\documentstyle{article}

\input amssym.def
\input amssym

\begin{document}

\title{Geometric zeta-functions on p-adic groups\footnote{AMS-class: 11F72, 11F75, 22E35, 22E40; Keywords: Selberg zeta function, continuous cohomology.}} 
\author{{\small by}\\{}\\  Anton Deitmar}
\date{}
\maketitle

\pagestyle{myheadings}
\markright{GEOMETRIC ZETA ON P-ADIC GROUPS}

\def \a{{{\frak a}}}
\def \ad{{\rm ad}}
\def \al{\alpha}
\def \ar{{\alpha_r}}
\def \A{{\Bbb A}}
\def \Ad{{\rm Ad}}
\def \b{{{\frak b}}}
\def \bs{\backslash}
\def \B{{\cal B}}
\def \cent{{\rm cent}}
\def \cov{{\rm cov}}
\def \C{{\Bbb C}}
\def \CA{{\cal A}}
\def \CB{{\cal B}}
\def \CD{{\cal D}}
\def \CE{{\cal E}}
\def \CF{{\cal F}}
\def \CG{{\cal G}}
\def \Chi{{\bf X}}
\def \CH{{\cal H}}
\def \CC{{\cal C}}
\def \CHC{{\cal HC}}
\def \CHS{{\cal HS}}
\def \CI{{\cal I}}
\def \CL{{\cal L}}
\def \CM{{\cal M}}
\def \CN{{\cal N}}
\def \CP{{\cal P}}
\def \CQ{{\cal Q}}
\def \CO{{\cal O}}
\def \CS{{\cal S}}
\def \CU{{\cal U}}
\def \det{{\rm det}}
\def \dist{{\rm dist}}
\def \e{\epsilon}
\def \End{{\rm End}}
\def \Fx{{\frak x}}
\def \FX{{\frak X}}
\def \g{{{\frak g}}}
\def \ga{\gamma}
\def \Ga{\Gamma}
\def \Gr{{\rm Gr}}
\def \h{{{\frak h}}}
\def \Hom{{\rm Hom}}
\def \im{{\rm im}}
\def \ind{{\rm ind}}
\def \Id{{\rm Id}}
\def \Im{{\rm Im}}
\def \Ind{{\rm Ind}}
\def \jtoinfty{\begin{array}{c} j\rightarrow \infty \\ \longrightarrow	\\ {}
		\end{array}}
\def \k{{{\frak k}}}
\def \K{{\cal K}}
\def \la{\lambda}
\def \lap{\triangle}
\def \m{{{\frak m}}}
\def \mod{{\rm mod}}
\def \n{{{\frak n}}}
\def \name{\bf}
\def \N{\Bbb N}
\def \o{{\frak o}}
\def \ord{{\rm ord}}
\def \O{{\cal O}}
\def \p{{{\frak p}}}
\def \ph{\varphi}
\def \prf{{\bf Proof: }}
\def \Pos{{\rm Pos}}
\def \qed{\hfill {$\Box$} 

$ $

}
\def \Q{\Bbb Q}
\def \res{{\rm res}}
\def \Res{{\rm Res}}
\def \R{{\Bbb R}}
\def \Re{{\rm Re \hspace{1pt}}}
\def \ra{\rightarrow}
\def \rank{{\rm rank}}
\def \supp{{\rm supp}}
\def \Si{\Sigma}
\def \Spec{{\rm Spec}}
\def \t{{{\frak t}}}
\def \T{{\Bbb T}}
\def \tr{{\hspace{1pt}\rm tr\hspace{1pt}}}
\def \vol{{\rm vol}}
\def \V{{\cal V}}
\def \z{{\frak z}}
\def \Z{\Bbb Z}
\def \={\ =\ }

\newcommand{\rez}[1]{\frac{1}{#1}}
\newcommand{\der}[1]{\frac{\partial}{\partial #1}}
\newcommand{\binom}[2]{\left( \begin{array}{c}#1\\#2\end{array}\right)}
\newcommand{\norm}[1]{\parallel #1 \parallel}

\newcounter{lemma}
\newcounter{corollary}
\newcounter{proposition}
\newcounter{theorem}

\renewcommand{\subsection}{\stepcounter{subsection}\stepcounter{lemma} 
	\stepcounter{corollary} \stepcounter{proposition}
	\stepcounter{conjecture}\stepcounter{theorem}
	$ $ \newline
	{\bf \arabic{section}.\arabic{subsection}\hspace{8pt}}}

%\vspace{4pt}

\newtheorem{conjecture}{\stepcounter{lemma} \stepcounter{corollary} 	
	\stepcounter{proposition}\stepcounter{theorem}
	\stepcounter{subsection}Conjecture}[section]
\newtheorem{lemma}{\stepcounter{conjecture}\stepcounter{corollary}	
	\stepcounter{proposition}\stepcounter{theorem}
	\stepcounter{subsection}Lemma}[section]
\newtheorem{corollary}{\stepcounter{conjecture}\stepcounter{lemma}
	\stepcounter{proposition}\stepcounter{theorem}
	\stepcounter{subsection}Corollary}[section]
\newtheorem{proposition}{\stepcounter{conjecture}\stepcounter{lemma}
	\stepcounter{corollary}\stepcounter{theorem}
	\stepcounter{subsection}Proposition}[section]
\newtheorem{theorem}{\stepcounter{conjecture} \stepcounter{lemma}
	\stepcounter{corollary}\stepcounter{proposition}		
	\stepcounter{subsection}Theorem}[section]
	
$$ $$

{\small
{\bf ABSTRACT.} We generalize the theory of $p$-adic geometric zeta functions  of Y. Ihara and K. Hashimoto to the higher rank case.
We give the proof of rationality of the zeta function and the connection of the divisor to group cohomology, i.e. the $p$-adic analogue of the Patterson conjecture.}

$$ $$

{\bf Introduction.} In \cite{ihara1} and \cite{ihara2} Y. Ihara defined geometric zeta functions for the group $PSL_2$.
This was the $p$-adic counterpart of the famous zeta functions of Selberg for Riemannian surfaces.
Using his results on the structure of discrete subgroups of $PSL_2$ of a $p$-adic field he proved, among other things, the rationality of these zeta functions.
By means of the geometric interpretation on the building,
K. Hashimoto \cite{hash1}, \cite{hash2}, \cite{hash3} extended this to arbitrary rank one groups.
In the present paper we generalize Ihara's theory to the higher rank case.
Our method uses harmonic analysis and is insofar not geometrical.
We will provide geometric interpretations in subsequent work.
We prove rationality of the zeta functions and show that the poles and zeroes are related to the cohomology of the discrete group, i.e. the $p$-adic counterpart of the Patterson-conjecture \cite{BuOl}.

\tableofcontents

\section{The trace formula}

\subsection
Let $F$ be a nondiscrete nonarchimedean locally compact field.
Then $F$ is either a finite extension of $\Q_p$ for some prime number $p$ or a field of formal power series in one indeterminate with coefficients in a finite field \cite{Weil}.
In any case $F$ is complete under a discrete valuation
$$
v : F^* \ra \Z .
$$

Let $\o$ be its ring of integers, $\m$ its maximal ideal and $\bar{F}=\o/\m$ its residue field.
Then $\bar{F}$ is finite, of characteristic, say, $p$ and has $q=p^n$ elements.
Define the absolute value $|.|$ on $F$ by $|x| := q^{-v(x)}$.
Note that this gives the modulus of $x$ with respect to any additive Haar-measure on $F$.

\subsection \label{haarmasse}
Let $\CG$ be a semi-simple linear algebraic group over ${F}$ and write $G$ for the group of $F$-rational points, so $G=\CG (F)$ is a totally disconnected unimodular locally compact group.
We will call an element $g$ of $G$ {\bf elliptic} if it lies in a compact subgroup of $G$.

We will fix Haar-measures of $G$ and its reductive subgroups as follows.
For $H\subset G$ being a torus there is a maximal compact subgroup $U_H$ which is open and unique up to conjugation. 
Then we fix a Haar measure on $H$ such that $\vol (U_H)=1$.
If $H$ is connected reductive with compact center then we choose the unique positive Haar-measure which up to sign coincides with the Euler-Poincar\'e measure
\cite{Kottwitz}.
So in the latter case our measure is determined by the following property:
For any discrete torsionfree cocompact subgroup $\Ga_H\subset H$ we have
$$
\vol (\Ga_H \bs H) \= (-1)^{q(H)}\chi(\Ga ,\Q),
$$
where $q(H)$ is the $k$-rank of the derived group $H_{der}$ and $\chi(\Ga ,\Q)$ the Euler-Poincar\'e characteristic of $H^.(\Ga ,\Q)$.
For the applications recall that centralizers of tori in connected groups are connected \cite{borel-lingroups}.

In the sequel we will fix a good maximal compact subgroup $K$ of $G$.

\subsection
Assume we are given a discrete subgroup $\Ga$ of $G$ such that the quotient space $\Ga \bs G$ is compact.
Let $(\omega ,V_\omega)$ be a finite dimensional unitary representation of $\Ga$ and let $L^2(\Ga \bs G,\omega)$ be the Hilbert space consisting of all measurable functions $f: G \ra V_\omega$ such that $f(\ga x) = \omega(\ga) f(x)$ and $|f|$ is square integrable over $\Ga \bs G$ (modulo null functions).
Let $R$ denote the unitary representation of $G$ on $L^2(\Ga \bs G,\omega)$ defined by right shifts, i.e. $R(g) \ph (x) = \ph (xg)$ for $\ph \in L^2(\Ga \bs G,\omega)$.
It is known that as a $G$-representation this space splits as a topological direct sum:
$$
L^2(\Ga \bs G,\omega) \= \bigoplus_{\pi \in \hat{G}} N_{\Ga ,\omega}(\pi) \pi
$$
with finite multiplicities $N_{\Ga ,\omega}(\pi)<\infty$.

We say that a function $f$ on $G$ is {\bf uniformly locally constant} if there exists a compact open subgroup $U$ of $G$ such that $f$ factors over $U\bs G/U$.

Let $f$ be integrable over $G$, so $f$ is in $L^1(G)$.
The integral
$$
R(f) := \int_G f(x) R(x) \ dx
$$
defines an operator on the Hilbert space $L^2(\Ga \bs G,\omega)$.

For $g\in G$ and $f$ any function on $G$ we define the {\bf orbital integral}
$$
\CO_g(f) := \int_{G_g\bs G} f(x^{-1} gx)\ dx,
$$
whenever the integral exists.
Here $G_g$ is the centralizer of $g$ in $G$.
In section \ref{haarmasse} we have fixed a Haar measure on $G_g$.
It is known that the group $G_g$ is unimodular, so we have an invariant measure on $G_g\bs G$.

\begin{proposition}
(Trace formula)
Let $f$ be integrable and uniformly locally constant, then we have 
$$
\sum_{\pi \in\hat{G}} N_{\Ga ,\omega}(\pi)\ \tr \pi (f) \= \sum_{[\ga]} \tr \omega (\ga)\ \vol(\Ga_\ga \bs G_\ga)\ \CO_\ga (f),
$$
where the sum on the right hand side runs over the set of $\Ga$-conjugacy classes $[\ga]$ in $\Ga$ and $\Ga_\ga$ denotes the centralizer of $\ga$ in $\Ga$.
Both sides converge absolutely and the left hand side actually is a finite sum.
\end{proposition}

\prf
We feel the necessity of giving a proof of this statement since, in contrast to the statement found in literature, we do not insist that $f$ has compact support. 

At first fix a fundamental domain $\CF$ for $\Ga \bs G$ and let $\ph\in L^2(\Ga \bs G,\omega)$, then
\begin{eqnarray*}
R(f) &\=& \int_Gf(y) \ph(xy)\ dy\\
	&\=& \int_G f(x^{-1}y) \ph(y)\ dy\\
	&\=& \sum_{\ga\in\Ga} \int_\CF f(x^{-1}\ga y)\ph(\ga y)\ dy\\
	&\=& \int_{\Ga\bs G} \left( \sum_{\ga\in\Ga} f(x^{-1}\ga y)\omega (\ga)\right) \ph (y) \ dy
\end{eqnarray*}

We want to show that the sum $\sum_{\ga\in\Ga} f(x^{-1}\ga y)\omega (\ga)$ converges in $\End(V_\omega)$ absolutely and uniformly in $x$ and $y$.
Since $y$ can be replaced by $\ga y$, $\ga\in\Ga$ and since $\omega$ is unitary, we only have to show the convergence of $\sum_{\ga\in\Ga}|f(x^{-1}\ga y)|$ locally uniformly in $y$.
Let $\ga$ and $\tau$ be in $\Ga$  and assume that $x^{-1}\ga y$ and $x^{-1}\tau y$ lie in the same class in $G/U$.
Then it follows $\tau yU\cap\ga yU \ne \emptyset$ so with $V=yUy^{-1}$ we have $\ga^{-1} \tau V \cap V \ne \emptyset$.
It is clear that $V$ depends on $y$ only up to $U$ so to show locally uniform convergence in $y$ it suffices to fix $V$.
Since $V$ is compact also $V^2 = \{ vv' |v,v'\in V\}$ is compact and so $\Ga \cap V^2$ is finite.
This implies that there are only finitely many $\ga\in\Ga$ with $\ga V\cap V\ne \emptyset$.
Hence the map $\Ga \ra G/U$, $\ga \mapsto x^{-1}\ga yU$ is finite to one with fibers having $\le n$ elements for some natural number $n$.
For $y$ fixed modulo $U$ we get
\begin{eqnarray*}
\sum_{\ga\in\Ga} |f(x^{-1} \ga y)| &\ \le\ & n\int_{G/U} |f(x)|\ dx\\
	&\=& \frac{n}{\vol(U)} \parallel f\parallel^1.
\end{eqnarray*}
We have shown the uniform convergence of the sum
$$
k_f(x,y) \= \sum_{\ga\in\Ga} f(x^{-1} \ga y) \omega(\ga).
$$

Observe that $R(f)$ factors over $L^2(\Ga\bs G,\omega)^U = L^2(\Ga\bs G/U,\omega)$, which is finite dimensional since $\Ga\bs G/U$ is a finite set.
So $R(f)$ acts on a finite dimensional space and $k_f(x,y)$ is the matrix of this operator.
We infer that $R(f)$ is of trace class, its trace equals
$$
\sum_{\pi\in\hat{G}} N_{\Ga ,\omega}(\pi)\ \tr\pi(f),
$$
and the sum is finite.
Further, since $k_f(x,y)$ is the matrix of $R(f)$ this trace also equals
\begin{eqnarray*}
\int_{\Ga \bs G} \tr k_f(x,x)\ dx &\=& \sum_{\ga\in\Ga} \int_\CF f(x^{-1}\ga x) \ dx\ \tr\omega(\ga)\\
	&\=& \sum_{[\ga]} \sum_{\sigma \in \Ga_\ga \bs \Ga} \int_\CF f((\sigma x)^{-1} \ga (\sigma x))\ dx\ \tr\omega(\ga)\\
	&\=& \sum_{[\ga]}  \int_{\Ga_\ga \bs G} f(x^{-1} \ga x)\ dx\ \tr\omega(\ga)\\
	&\=& \sum_{[\ga]}  \vol(\Ga_\ga\bs G_\ga)\ \int_{G_\ga \bs G} f(x^{-1} \ga x)\ dx\ \tr\omega(\ga)
\end{eqnarray*}
by formal computation.
\qed

\section{Special test functions}
Fix a $k$-parabolic of rank one, i.e. $P=LN$, where $N$ is the unipotent radical of $P$ and $L$ a Levi-component, such that
the maximal split torus $A$ in the center of $L$ is  one dimensional.
Note that such $P$ always exist unless $G$ is compact.
We can write $P=MAN$, $M$ is reductive with compact center and the centralizer of $A$ equals $L=MA$.
Let $\bar{P}=MA\bar{N}$ be the opposite parabolic and write $\n$ resp. $\bar{\n}$ for the Lie algebras of $N$ and $\bar{N}$.

Let $A^+$ be the subset of those $a\in A$ such that $a$ acts on the Lie algebra of $N$ by eigenvalues $\mu$ such that $|\mu |>1$.
This fixes an isomorphism $A\ra F^*$ such that $A^+$ is mapped onto the elements of absolute value $>1$.
Let accordingly $A^-$ denote the preimage of $\m\cap k^*$.
It may happen that $A^+$ and $A^-$ are conjugate in $G$.
Via this isomorphism $A\ra F^*$ we define the absolute value $|.|:A\ra \R_+^*$ on $A$.
This in turn defines a valuation on $A^*=\Hom_{cont}(A,\C^*)$, the group of quasi-characters of $A$, as follows: let $\chi$ be in $A^*$, taking absolute value in $\C$ we get a character from $A$ to $\R_+^*$ of the form
$$
|\chi (a)| = |a|^c,
$$
for some $c\in\R$.  
Define
$$
|\chi |\ :=\ |c|.
$$
Any root $\alpha$ of $A$ will be viewed as a homomorphism $A\ra F^*$, $a\mapsto a^\alpha$. 
Taking absolute values induces a homomorphism $A\ra \R_+^*$ denoted as $a\mapsto |a|^{c(\alpha)}$.
Define $|\alpha |:= |c(\alpha )|$.

Let $G$ act on itself by conjugation, write $g.x = gxg^{-1}$, write $G.x$ for the orbit, so $G.x = \{ gxg^{-1} | g\in G \}$ as well as $G.S = \{ gsg^{-1} | s\in S , g\in G \}$ for any subset $S$ of $G$.
We are going to consider functions that are supported on the closure of the set $G.(MA)$. Now let $f_{EP}$ be an Euler-Poincar\'e function to the group $M$.
We will not choose $f_{EP}$ as in \cite{Kottwitz} since we need $f_{EP}$ to be invariant under $K_M=K\cap M$-conjugation.
So let $\tilde{f}_{EP}$ be the Euler-Poincar\'e function of Kottwitz then we define
$$
f_{EP}(m) \ \:= \ \rez{\vol(K_M)} \int_{K_M} \tilde{f}_{EP} (kmk^{-1}) dk.
$$
It follows that $f_{EP}$ still is locally constant of compact support and that it has the same orbital integrals and traces as $\tilde{f}_{EP}$.

\subsection
For $g\in G$ and $V$ any $F$ vector space on which $g$ acts linearly let $E(g|V)\subset \R^*_+$ be defined by
$$
E(g|V) \ :=\ \{ |\mu | : \mu {\rm \ is\ an\ eigenvalue\ of}\ g {\rm \ on\ } V\}.
$$
Here we have extended the absolute value $|.|$ to the algebraic closure $\bar{F}$ of $F$.
Note that, since the $g$-action on $V$ is defined over $F$, the set $E(g|V)$ does not depend on the choice of the extension.
Let $\la_{min}(g|V):= \min(E(g|V))$ and  $\la_{max}(g):= \max(E(g|V))$ the minimum and maximum.
Define
$$
\la(am) := \frac{\la_{min}(a|\n)}{\la_{max}(m|\g)},
$$
where $\g$ and $\n$ are the Lie algebras of $G$ and $N$ adn $\alpha$ is the short positive root in $\Phi(\g ,\a)$.
Note that $\la(am)$ is well defined although it may happen that $am=a'm'$ with $a\ne a'$, because then $a$ and $a'$ only differ by an elliptic element.
We will construct a function on the set
$$
(AM)^{\sim} \ :=\ \{ am\in AM | \la(am)>1 \}.
$$
The following properties of $(AM)^{\sim}$ are immediate
\begin{itemize}
\item[1.]
$A^+M_{ell}\subset (MA)^{\sim}$
\item[2.]
$am\in (AM)^{\sim} \Rightarrow a\in A^+$
\item[3.]
$am, a'm' \in (AM)^{\sim}, g\in G\ {\rm with}\ a'm'=gamg^{-1} \Rightarrow a=a', g\in AM$.
\end{itemize}
Here we have written $M_{ell}$ for the set of elliptic elements in $M$.

For $T\in \C$ consider  the function $g_T$ on $(AM)^{\sim}$ given by
$g_T(am) = T^{l_{am}}$, where $l_{am}$ is defined by
$$
l_{am} \ :=\ \frac{\log (\la (am))}{|\alpha |\log(q)}.
$$
Here $\alpha$ is the short positive root in the root system $\Phi(\a ,\g)$.
Note that $l_{am}$ is an integer and
that for $m$ elliptic we have $l_{am}=l_a=v(a)$, the {\bf length} of $a$.
Choose any locally constant $\eta : N \ra \R$ which has compact support, is positive, invariant under $K\cap M$ and such that
$\int_N \eta(n) dn =1$.

Fix a finite dimensional complex representation $(\sigma ,V_\sigma)$ of $M$,
further fix a unitary character $\chi$ of $A$.

Given these data, let $\Phi =\Phi_{\eta, \sigma ,\chi ,T} : G\ra \C$ be defined by
$$
\Phi(kn ma (kn)^{-1}) \= \eta(n) f_{EP}(m) \tr \sigma(m)\chi(a){g_T(am)},
$$
for $k\in K$, $n\in N$, $am\in (AM)^{\sim}$. Further $\omega(g)=0$ if $g$ is not in $G.(MA)^{\sim}$.

In order to see that $\Phi$ is well-defined recall first that by the decomposition $G=KP=KNMA$ every $g\in G.(MA)^{\sim}$ can be written in the form $kn ma (kn)^{-1}$.
By the properties of $(AM)^{\sim}$ we see that two those representations can only differ by an element of $K\cap P =K\cap M$.

The invariance of $\eta$ and of $f_\sigma$ under $K\cap M$-conjugation shows that $\Phi$ is well-defined with respect to the $K$-conjugation.

Let $\CE_P(\Ga)$ be the set of all $\Ga$-conjugacy classes $[\ga ]$ such that $\ga$ is in $G$ conjugate to some element $a_\ga m_\ga$ of $A^+M_{ell}$.

\begin{proposition} \label{goes_in_tf}
The function $\Phi$ is uniformly locally constant and lies in $L^1(G)$ for  $|T|$ sufficiently small. 
For any finite dimensional unitary representation $\omega$ of $\Ga$ we have
$$
\sum_{\pi \in \hat{G}} N_{\Ga,\omega} (\pi)\ \tr\pi(\Phi) \= \sum_{[\ga]\in\CE_P(\Ga)} \vol (\Ga_\ga \bs G_\ga )\ 
{\tr (\omega(\ga))\tr\sigma(m_\ga)\chi(a_\ga)\ T^{l_\ga}},
$$
where the sum runs over all classes $[\ga]$ such that $\ga$ is conjugate to an element $m_\ga a_\ga$ of $M_{ell}A^+$ and we have written $l_\ga$ for the length $l_{a_\ga}$.
\end{proposition}

\prf
We recall the Weyl integration formula:
Let $H_1,\dots ,H_r$ denote a maximal set of  Cartan-subgroups of $G$ which are pairwise non conjugate.
For any $j$ let $N_G(H_j)$ be the normalizer of $H_j$ in $G$ and $Z_G(H_j)$ the centralizer.
Let $W(G,H_j):= N_G(H_j)/Z_G(H_j)$ be the Weyl group, which is finite.
For $h\in H_j$ let
$$
D_j(h)\ :=\ \left| \det(1-h |\g /\h_j)\right| ,
$$
the Jacobian, where $\h_j$ is the Lie algebra of $H_j$.
The Weyl integration formula states that for $f\in L^1(G)$ the integral $\int_G f(x)\ dx$ equals
$$
\sum_{j=1}^r \rez{|W(G,H_j)|} \int_{H_j} \int_{G/H_j} f(ghg^{-1}) D_j(h)\ dgdh.
$$
It also shows that if $f$ is measurable and the right hand side exists, then $f$ is integrable.

Recall that an element $g$ of $G$ is called {\bf regular} if its centralizer is a torus.
Note that the set $G^{reg}$ of regular elements in $G$ is open and dense and that the Weyl integration formula remains true if we restrict the integration domains to the subsets of regular elements.

For any subset $S$ of $G$ set $S^{reg}:=S\cap G^{reg}$.
Assume that our Cartan subgroups $H_1,\dots ,H_r$ are ordered in a way that
$H_1^{reg},\dots ,H_s^{reg}$ meet $G.(MA)$ and $H_{s+1}^{reg},\dots ,H_r^{reg}$ don't.
Replacing the $H_j$'s by conjugates if necessary we assume $H_j\subset MA$ for $j\le s$.
This implies $H_j = H_j^M A$ for $j\le s$, where the $H_j^M$ run through a set of nonconjugate Cartan subgroups of $M$.
This set needn't be maximal since two Cartans $H,H'\subset MA$ could be conjugate in $G$ without being conjugate in $MA$.
There is, however, a unique $MA$-conjugacy class of Cartans $H=BA$ with $B\subset M$ compact.
So we may assume $H_1=BA$ and $H_j^M$ noncompact for $j>1$.
Note that this implies that any $h\in H_j^{M,reg}$ is nonelliptic, where regular here means regular in $M$.

Since our test function $\Phi$ is supported in $G.(MA)$ it follows that for $j>s$ we have
$$
\int_{H_j}\int_{G/H_j} \Phi(ghg^{-1})\ dgdh \= 0.
$$
For $j=2,\dots ,s$ we have
$$
\int_{H_j}\int_{G/H_j} \Phi(ghg^{-1})D_j(h)\ dgdh
$$
\begin{eqnarray*}
	&\=& \int_{H_j}\int_{KNM/H_j^M} \Phi((knm)h(knm)^{-1}) D_j(h)\ dkdndmdh\\
	&\=& \int_{H_j} \CO_{h^M}^M(f_{EP})\ 
		{g_T(h)} \tr\sigma(h^M)\chi(h^A)D_j(h) \ dh,
\end{eqnarray*}
where $\CO_{h^M}^M(f)$ denotes the $M$-orbital integral of $h^M$, the projection to $M$ and $h^A$ is the projection to $A$.
Since $h^M$ is not elliptic for $h\in H_j^{reg}$ it follows that this term vanishes \cite{Kottwitz}.

It remains to consider $H_1 =H=AB$.
For the term corresponding to $H$ we have that
$$
\int_{H_1}\int_{G/H_1} \Phi(ghg^{-1}) D_j(h) \ dgdh
$$
equals
$$
|W(G,A)| \int_{A^+B} {T^{l_a}}\tr\sigma(b)
\chi(a)D_j(ab)\ dadb.
$$
It is easy to see that this integral exists for $|T|$ small enough.
Therefore we conclude that $\Phi\in L^1(G)$.
Since it clearly is uniformly locally constant, the trace formula applies to give:
$$
\sum_{\pi \in\hat{G}} N_{\Ga ,\omega}(\pi)\ \tr \pi (\Phi) \= \sum_{[\ga]} \tr \omega (\ga)\ \vol(\Ga_\ga \bs G_\ga)\ \CO_\ga (\Phi).
$$
In the above argumentation we in fact showed that $\CO_\ga(\Phi)=0$ unless $\ga$ is $G$-conjugate to an element $m_\ga a_\ga \in M_{ell}A^+$.
In the latter case we compute
$$
\CO_\ga (\Phi) \= \CO_{m_\ga a_\ga}(\Phi)\= \CO_{m_\ga}^M(f_{EP})\tr\sigma(m_\ga)\chi(a_\ga)
{T^{l_\ga}}.
$$
From \cite{Kottwitz} we now take that the orbital integrals $\CO_{m_\ga}^M(f_{EP})$ are all equal to $1$.
\qed

\section{The covolume of a centralizer}

Suppose $\ga\in\Ga$ is $G$-conjugate to some $a_\ga m_\ga\in A^+M_{ell}$.
We want to compute the covolume
$$
\vol(\Ga_\ga \bs G_\ga) .
$$
We suppose that $\Ga$ is neat. This implies that for any $\ga\in\Ga$ the Zariski closure of the group generated by $\ga$ is a torus.
It then follows  that $G_\ga$ is a connected reductive group \cite{borel-lingroups}.

An element $\ga\in\Ga$ is called {\bf primitive} if $\ga =\sigma^n$ with $\sigma\in\Ga$ and $n\in\N$ implies $n=1$.
It is a property of discrete cocompact torsion free subgroups $\Ga$ of $G$ that every $\ga\in\Ga$, $\ga\ne 1$ is a positive power of a unique primitive element.
In other words, given a nontrivial $\ga\in\Ga$ there exists a unique primitive $\ga_0$ and a unique $\mu(\ga)\in\N$ such that
$$
\ga =\ga_0^{\mu(\ga)}.
$$

Let $\Sigma$ be a group of finite cohomological dimension over $\Q$.
We write
$$
\chi(\Sigma) \= \chi(\Sigma ,\Q) \ :=\ \sum_{p=0}^\infty (-1)^p \dim H^p(\Sigma ,\Q),
$$
for the Euler-Poincar\'e characteristic.
In \cite{onsome} we established the notion of higher Euler characteristics $\chi_{_n}$ for $n\in \N$.
The first one of these is
$$
\chi_{_1}(\Sigma) \= \chi_{_1}(\Sigma ,\Q) \ :=\ \sum_{p=0}^\infty p(-1)^{p+1} \dim H^p(\Sigma ,\Q).
$$

We call the group $\Ga$ {\bf neat} if for every $\ga\in\Ga$ and every rational representation $\rho : G\ra GL_n(F)$ of $G$ the matrix $\rho(\ga)$ has no eigenvalue which is a root of unity.
Every arithmetic $\Ga$ has a finite index subgroup which is neat \cite{borel}.

\begin{proposition}
Assume $\Ga$ neat and let $\ga\in\Ga$ be $G$-conjugate to an element of $A^+M_{ell}$.
Then we get
$$
\vol (\Ga_\ga\bs G_\ga) = l_{\ga_0}\ (-1)^{q(G)+1}\ \chi_{_1}(\Ga_\ga).
$$
\end{proposition}

\prf
Modulo $G$-conjugation we may assume $\ga = a_\ga m_\ga\in A^+M_{ell}\subset (AM)^\sim$.
By property 3. of $(AM)^\sim$ we conclude
$$
G_\ga =AM_{m_\ga}.
$$
Since $m_\ga$ is elliptic we conclude that $M_{m_\ga}$ has no central split torus.
So we can write $M_{m_\ga}= CM_{m_\ga}'$, where $C$ is compact abelian and $M_{m_\ga}'$ is semi-simple.
We have $\ga\in AC$ and $\Ga \cap (AC)$ is generated by $\ga_0$.
Our measure normalizations imply that $\vol(C)=1$ and we get
\begin{eqnarray*}
\vol (\Ga_\ga\bs G_\ga) &\=& \vol(C\Ga_\ga \bs G_\ga)\\
	&\=& \vol (U_A C \Ga_\ga \bs ACM_{m_\ga}'),
\end{eqnarray*}
where $U_A$ is the maximal compact subgroup of $A$.
Let $\sigma\in\Ga_\ga$ such that $\sigma = a_\sigma m_\sigma$ with $a_\sigma \in A^+$ having minimal absolute value.
Then $|a_\sigma|^k = |a_{\ga_0}|$ for some natural number $k$
It follows $a_\sigma^k = \ga_0$ modulo $U_A C$.
So let $\Ga_\ga^1$ be the group generated by $\Ga_\ga$ and $a_\sigma$.
We find that
$$
[\Ga_\ga^1 U_A C : \Ga_\ga U_A C] \= k.
$$
Let further $\Sigma := M_{m_\ga}'\cap AC\Ga_\ga$. 
Then $\Sigma$ is discrete and cocompact in $M_{m_\ga}'$.
Since $\Ga$ is neat it is also torsion free modulo the center $Z$ of $M_{m_\ga}'$ which is finite and contained in $\Sigma$.
It follows
\begin{eqnarray*}
\vol (\Ga_\ga \bs G_\ga) &\=& \vol (U_A C \Ga_\ga^1 \bs AC M_{m_\ga}')\\
	&\=& k\ \vol(U_A C \langle a_\sigma \rangle \bs AC \ \times\ \Sigma \bs M_{m_\ga}')\\
	&\=& k\ l_{a_\sigma} \vol(\Sigma \bs M_{m_\ga}')\\
	&\=& l_{\ga_0} (-1)^{q(M_{m_\ga}')}\ \chi (\Sigma /Z)
\end{eqnarray*}
Now $q(G) = q(M_{m_\ga}')+1$ and it remains to show that
$$
\chi (\Sigma /Z) \= \chi_{_1}(\Ga_\ga).
$$

Since $\Ga_\ga \cap AC = \langle \ga_0 \rangle$ it follows that we have an exact sequence of groups
$$
1\ra \langle \ga_0 \rangle \ra \Ga_\ga \ra \Sigma /Z \ra 1.
$$
Since $\Ga_\ga$ and $\Sigma /Z$ are cocompact torsion free subgroups of reductive $p$-adic groups, they are of finite cohomological dimension over $\Q$.
Our proposition follows from the 

\begin{lemma}\label{chichi1}
Let $G,Q$ be of finite cohomological dimension over $\Q$.
Let $C$ be an infinite cyclic group and assume there is an exact sequence
$$
1\ra C\ra G\ra Q\ra 1.
$$
Then we have
$$
\chi (Q,\Q) \= \chi_{_1}(G,\Q).
$$
\end{lemma}

\prf
Consider the Hochschild-Serre spectral sequence:
$$
E_2^{p,q} \= H^p(Q,H^q(C,\Q))
$$
which abuts to
$$
H^{p+q}(G,\Q).
$$
Since $C\cong \Z$ it follows
$$
H^q(C,\Q) \= \left\{ \begin{array}{cl} \Q& {\rm if}\ q=0,1\\ 0&{\rm else}.\end{array}\right.
$$
So the spectral sequence degenerates and the claim is clear.
\qed

\section{Rationality of the zeta function}

We keep fixed a neat cocompact discrete subgroup $\Ga$ of $G$ and a rank one $F$-parabolic $P=MAN$.
Write $\CE_P^p(\Ga)$ for the set of all $[\ga]\in\CE_P(\Ga)$ such that $\ga$ is primitive.
For a finite dimensional unitary representation $(\omega ,V_\omega)$ of $\Ga$ consider the infinite product
$$
Z_{P,\sigma ,\chi ,\omega} (T) \ :=\ \prod_{[\ga ]\in\CE_P^p(\Ga)} \det\left( 1-T^{l_\ga}\chi(a_\ga)\omega (\ga)\otimes\sigma(m_\ga)\right)^{\chi_{_1}(\Ga_\ga)}.
$$

\begin{theorem}
The infinite product $Z_{P,\sigma ,\chi ,\omega}(s)$ converges for $|T|$ small enough. 
The limit extends to a rational function in $T$.
\end{theorem}

\prf
Formally at first we compute the logarithmic derivative as
$$
\frac{Z_{P,\sigma ,\chi ,\omega}'}{Z_{P,\sigma ,\chi ,\omega}}(T) \= \frac{d}{dT} \left(-\sum_{[\ga]\in\CE_P^p(\Ga)} \chi_{_1}(\Ga_\ga) \sum_{n=1}^\infty \frac{T^{nl_\ga}}{n} \tr \omega(\ga^n)\tr\sigma(m_\ga^n)\chi(a_\ga^n)\right)
$$ $$
= -\sum_{[\ga]\in\CE_P^p(\Ga)}\chi_{_1}(\Ga_\ga) \sum_{n=1}^\infty l_\ga \frac{T^{nl_\ga}}{T} \tr \omega(\ga^n)\tr\sigma(m_\ga^n)\chi(a_\ga^n).
$$
Since $\Ga$ is neat we have $\Ga_\ga =\Ga_{\ga^n}$ for all $n\in\N$, so
$$
\frac{Z_{P,\sigma ,\chi ,\omega}'}{Z_{P,\sigma ,\chi ,\omega}}(T) \= -\sum_{[\ga]\in\CE_P(\Ga)}\chi_{_1}(\Ga_\ga) l_{\ga_0} \frac{T^{l_\ga}}{T} \tr \omega(\ga)\tr\sigma(m_\ga)\chi(a_\ga)
$$ $$
= (-1)^{q(G)}\sum_{[\ga]\in\CE_P(\Ga)} \vol(\Ga_\ga\bs G_\ga)\frac{T^{l_\ga}}{T} \tr \omega(\ga)\tr\sigma(m_\ga)\chi(a_\ga)
$$
by the result of the previous section.
This last expression equals the right hand side of the trace formula of Proposition \ref{goes_in_tf} and so it converges for $\Re(s)>>0$.
A fortiori the convergence of the product follows.

We have proven
$$
\frac{Z_{P,\sigma ,\chi ,\omega}'}{Z_{P,\sigma ,\chi ,\omega}}(T) \= (-1)^{q(G)} \sum_{\pi\in\hat{G}} N_{\Ga ,\omega} (\pi)\ \frac{\tr\pi(\Phi)}{T}.
$$
The right hand side is a finite sum.
A result of Harish-Chandra says that for $\pi\in\hat{G}$ there is a locally integrable function $\Theta_\pi$ on $G$ such that for all locally constant functions $f$ of compact support we have
$$
\tr\pi(f) \= \int_G f(x) \Theta_\pi(x)\ dx.
$$
This function $\Theta_\pi$ is called the {\bf character} of $\pi$.
Recall that a vector $v\in V_\pi$, the representation space of $\pi$ is called {\bf smooth}
if its stabilizer $G_v$ in $G$ is open.
The space of smooth vectors $V_\pi^\infty$ is dense in $V_\pi$ and forms an {\bf admissible} representation of $G$, i.e. for any open subgroup $H$ of $G$ the space of fixed vectors $V_\pi^H$ is finite dimensional.
Let $P=MAN$ be our fixed parabolic, $\bar{P}=MA\bar{N}$ its opposite and write $\pi_{\bar{N}}$ for the {\bf Jacquet module} of $\pi$.
By definition, $\pi_{\bar{N}}$ is the largest quotient of $V_\pi^\infty$ on which $\bar{N}$ acts trivially.
One can achieve this by factoring out the submodule generated by all vectors of the form $v-\pi(n)v$, $v\in V_\pi^\infty$, $n\in N$.
It is known that $\pi_{\bar{N}}$ is a finitely generated admissible module for the 
group $MA$ and therefore it has a character $\Theta_{\pi_{\bar{N}}}^{MA}$.
In \cite{Casselman} it is shown that
$$
\Theta_\pi(am) \= \Theta_{\pi_{\bar{N}}}^{MA}(ma)
$$
for $ma\in A^+M_{ell}$.

Since the character $\Theta_\pi$ is conjugation invariant, i.e. $\Theta_\pi(gxg^{-1})=\Theta_\pi(x)$, the Weyl integration formula gives us
$$
\tr \pi (\Phi) \= \sum_{j=1}^r \rez{|W(G,H_j)|} \int_{H_j} \Theta_\pi(h) D_j(h) \int_{G/H_j} \Phi(ghg^{-1}) \ dg dh.
$$
As in the proof of Proposition \ref{goes_in_tf} we find that only the term with $H_1=AB$ will survive.
So
$$
\tr \pi(\Phi) \= \rez{|W(G,AB)|} \int_{AB} \Theta_\pi(ab) D_1(ab) \int_{G/AB} \Phi(gabg^{-1})\ dgdadb.
$$
First note that $\Theta_\pi(ab) = \Theta_{\pi_{\bar{N}}}^{MA}(ab)$ and then observe that the integral over $A$ may be splitted into an integral over $A^+$ and one over $A^-$, since $\Phi$ vanishes on the rest.
If $|W(G,A)|=2$ then these two coincide, else we have $|W(G,A)|=1$ and then $\Phi$ vanishes on $A^-$.
In any case we get
$$
\tr\pi(\Phi) \= \rez{|W(G,B)|} \int_{A^+}\int_B \Theta_{\pi_{\bar{N}}}^{MA}(ab) D_1(ab) \int_{G/AB} \Phi(gabg^{-1})\ dgdadb.
$$

We denote the Lie algebra of $N$ by $\n$ and the Lie algebra of $\bar{N}$ by $\bar{\n}$.
According to the definition of $\Phi$ the inner integral over $G/AB$ gives
$$
\CO_b^M(f_{EP}) g_T(ab) \= \CO_b^M(f_{EP}) g_T(a)\tr\sigma(b)\chi(a),
$$ 
where $\CO_b^M$ is the orbital integral with respect to the group $M$.
Furthermore we have
$$
D_1(ab) \= |\det(1-ab| \n\oplus \bar{\n})|.
$$
Now assume $a\in A^+$ then any eigenvalue of $ab$ on $\n$ is of valuation $>0$ and any eigenvalue on $\bar{\n}$ is of valuation $<0$.
By the ultrametric property of the valuation it follows:
\begin{eqnarray*}
D_1(ab) &=& |\det(1-ab|\n )|\\
	&=& |\det(ab|\n)|\\
	&=& |\det(a|\n)|\\
	&=& |a|^{2|\rho|}.
\end{eqnarray*}
Here we have written $\rho$ for the usual modular shift, i.e. the half sum of the positive roots.
Since the nonelliptic orbital integrals of the function $f_{EP}$ all vanish \cite{Kottwitz}, they can formally be added to what we have.
An application of the Weyl integration formula of $M$ yields
\begin{eqnarray*}
\tr\pi(\Phi) &=& \int_{A^+} \int_{M} \Theta_{\pi_{\bar{N}}}^{MA}(am) f_{EP}(m) \tr\sigma(m)\chi(a) |a|^{2|\rho|} g_T(a)\ dadm\\
	&=& \int_{A^+} \int_{M} \Theta_{\pi_{\bar{N}}\otimes\sigma_\chi}^{MA}(am) f_{EP}(m)  |a|^{2|\rho|} g_T(a)\ dadm,
\end{eqnarray*}
where we have written $\sigma_\chi$ for the representation of $AM$ given by $\sigma_\chi(am) = \chi(a)\sigma(m)$.

We write $H_d^i(M,\pi_{\bar{N}})$ for the differentiable cohomology groups (\cite{Borel-Wallach}, chap. X).
Casselman showed (\cite{Kottwitz}) that
$$
\tr \pi_{\bar{N}}\otimes\sigma_\chi (f_{EP}) \= \sum_i (-1)^i \dim H_d^i(M,\pi_{\bar{N}}\otimes \sigma_\chi).
$$
The cohomology groups $H_c^i(M,\pi_{\bar{N}})$ are finite dimensional $A$-modules and we get
$$
\tr \pi(\Phi) \= \sum_{i=0}^{\dim M} (-1)^i \int_{A^+} \tr (a|H_c^i(M,\pi_{\bar{N}}\otimes\sigma_\chi))\ 
|a|^{2|\rho|} T^{v(a)} \ da.
$$
Decompose the $A$-module $H_c^i(M,\pi_{\bar{N}}\otimes\sigma_\chi)$ into generalized eigenspaces under $A$.
Clearly only that part of $H_c^i(M,\pi_{\bar{N}}\otimes\sigma_\chi)$ survives on which $A$ acts by unramified characters.
Suppose this space splits as
$$
\bigoplus_{\la\in\C} E_\la^i,
$$
where for $\la\in\C$ the space $E_\la^i$ is the generalized eigenspace to the character $|a|^\la$.
Let $m_\la^i := \dim_\C E_\la^i$.
Then
\begin{eqnarray*}
\frac{\tr \pi(\Phi)}{T} &\= & \sum_{i=1}^{\dim M} (-1)^i \sum_{\la\in\C} m_\la^i \int_{A^+} |a|^{2|\rho|+\la} T^{v(a)-1}\ da\\
	&\=& \sum_{i=1}^{\dim M} (-1)^i \sum_{\la\in\C} m_\la^i \sum_{n=1}^\infty (Tq^{2|\rho|+\la})^n\ T^{-1}\\
	&\=& -\sum_{i=1}^{\dim M} (-1)^i \sum_{\la\in\C} m_\la^i  \rez{T-q^{-2|\rho|-\la}}.
\end{eqnarray*}
Since only finitely many $\pi\in\hat{G}$ contribute, it follows that $Z_{P,\sigma ,\chi ,\omega}'/Z_{P,\sigma ,\chi ,\omega}$ is a finite sum of functions of the form $\rez{T-a}$, this implies that  $Z_{P,\sigma ,\chi ,\omega}$ is a rational function. 
The theorem follows.
\qed

The proof of the theorem actually gives us:

\begin{proposition}
For $\la\in\C$ let $H_c^i(M,\pi_{\bar{N}}\otimes\sigma_\chi)_\la$ denote the generalized $|.|^\la$-eigenspace of the $A$-action on $H_c^i(M,\pi_{\bar{N}}\otimes\sigma_\chi)$, then the (vanishing-) order of $Z_{P,\sigma ,\chi ,\omega}(q^{-s-2|\rho |})$ in $s=s_0$  equals
$$
(-1)^{q(G)+1}\sum_{\pi\in\hat{G}} N_{\Ga ,\omega}(\pi) \sum_{i=0}^{\dim M} (-1)^i \dim H_d^i(M,\pi_{\bar{N}}\otimes\sigma_\chi)_{-s_0}.
$$
\end{proposition}
\qed

\section{The divisor of the zeta function and group cohomology}

First we introduce some standard notation.
A representation $(\rho ,V)$ on a complex vector space will be called {\bf algebraic} if every vector $v\in V$ is fixed by some compact open subgroup.
Let $Alg(G)$ be the category of algebraic representations (comp. \cite{BerZel}).

The representation $(\rho ,V)$ is called {\bf admissible} if for any compact open subgroup $K\subset G$ the space of fixed vectors $V^K$ is finite dimensional.
Let $\CHC_{alg}(G)$ denote the full subcategory of $Alg(G)$ consisting of all admissible representations of finite length in $Alg(G)$.

In the following we will mix the representation theoretic notion with the module-theoretic notion, so $V$ will be a $G$-module with $x\in G$ acting as $v\mapsto x.v$.

\subsection
We introduce the setup of topological (continuous) representations.
Let $\CC(G)$ denote the category of continuous representations of $G$ on locally convex, complete, Hausdorff topological vector spaces.
The morphisms in $\CC(G)$ are continuous linear $G$-maps.

To any $V\in \CC(G)$ we can form $V^\infty$, the algebraic part, which is by definition the space of all vectors in $V$ which have an open stabilizer.
By the continuity of the representation it follows that $V^\infty$ is dense in $V$.
The association $V\mapsto V^\infty$ gives an exact functor from $\CC(G)$ to $Alg(G)$.
An object $V\in \CC(G)$ is called admissible if $V^\infty$ is and so we define $\CHC_{top}$ to be the full subcategory of $\CC(G)$ consisting of all admissible representations of finite length.
Let $F$ be the functor:
\begin{eqnarray*}
F : \CHC_{top} &\ra & \CHC_{alg}\\
W&\mapsto & W^\infty .
\end{eqnarray*}

Let $V$ be in $\CHC_{alg}$ then any $W\in\CHC_{top}$ with $FW=V$ will be called a {\bf completion} of $V$.

To $V\in Alg(G)$ let $V^*=Hom_\C(V,\C)$ be the dual module and let $\tilde{V}$ be the algebraic part of $V^*$.
If $V$ is admissible we have that the natural map from $V$ to $\tilde{\tilde{V}}$ is an isomorphism.

On the space $C(G)$ of continuous maps from $G$ to the complex numbers we have two actions of $G$, the left and the right action given by
$$
L(g)\ph(x):=\ph(g^{-1}x)\ \ \  R(g)\ph(x) := \ph(xg),
$$
where $g,x\in G$ and $\ph\in C(G)$. 
For $V\in\CHC_{alg}(G)$ let
$$
V^{max} := Hom_G(\tilde{V} ,C(G)),
$$
where we take $G$-homomorphisms with respect to the right action, so $V^{max}$ is the space of all linear maps $f$ from $\tilde{V}$ to $C(G)$ such that $f(g.v^*)=R(g)f(g^*)$.
Then $V^{max}$ becomes a $G$-module via the left translation: for $\alpha\in V^{max}$ we define $g.\alpha (v^*) := L(g)\alpha(v^*)$.
We call $V^{max}$ the {\bf maximal completion} of $V$.
The next lemma and the next Proposition will justify this terminology.
First observe that the topology of locally uniform convergence may be installed on $V^{max}$ to make it an element of $\CC(G)$.
By $V\mapsto V^{max}$ we then get a functor $R:\CHC_{alg}(G)\ra\CHC_{top}(G)$.

\begin{lemma}
We have $(V^{max})^\infty \cong V$.
\end{lemma}

\prf
There is a natural map $\Phi : V\ra V^{max}$ given by $\Phi(v)(v^*)(g) := v^*(g^{-1}.v)$, for $v^*\in \tilde{V}$ and $g\in G$.
This map is clearly injective.
To check surjectivity let $\alpha\in(V^{max})^\infty$, then the map
$\tilde{V}\ra \C$, $v^*\mapsto \alpha(v^*)(1)$ lies in $\tilde{\tilde{V}}$.
Since the natural map $V\ra \tilde{\tilde{V}}$ is an isomorphism, there is a $v\in V$ such that $\alpha(v^*)(1)=v^*(v)$ for any $v^*\in \tilde{V}$, hence $\alpha(v^*)(g)= \alpha(g.v^*)(1)=g.v^*(v)=v^*(g^{-1}.v)=\Phi(v)(v^*)(g)$, hence $\alpha =\Phi(v)$.
\qed

It follows that the functor $R$ maps $\CHC_{alg}(G)$ to $\CHC_{top}(G)$ and that $FR=Id$.

\begin{proposition}
The functor $R: \CHC_{alg}(G)\ra \CHC_{top}(G)$ is right adjoint to $F: W\mapsto W^\infty$.
So for $W\in \CHC_{top}(G)$ and $V\in \CHC_{alg}(G)$ we have a functorial isomorphism:
$$
Hom_{alg}(FW,V)\cong Hom_{top}(W,RV).
$$
\end{proposition}

\prf
We have a natural map
\begin{eqnarray*}
Hom_{top}(W,V^{max}) &\ra& Hom_{alg}(W^\infty ,V)\\
\alpha &\mapsto& \alpha |_{W^\infty}.
\end{eqnarray*}
Since $W^\infty$ is dense in $W$ this map is injective.
For surjectivity we first construct a map $\psi : W\ra (W^\infty)^{max}$.
For this let $W'$ be the topological dual of $W$ and for $w'\in W'$ and $w\in W$ let
$$
\psi_{w',w}(x) := w'(x^{-1}.w),\ \ \ x\in G.
$$
The map $w\mapsto \psi_{.,w}$ gives an injection
$$
W\hookrightarrow Hom_G(W',C(G)).
$$
Let $\ph\in \tilde{FW}$, then $\ph$ factors over $(FW)^K$ for some compact open subgroup $K\subset G$.
But since $(FW)^K =W^K$ it follows that $\ph$ extends to $W$.
The admissibility implies that $\ph$ is continuous there.
So we get $\tilde{(FW)} \hookrightarrow W'$.
The restriction then gives
$$
Hom_G(W',C(G))\ra Hom_G(\tilde{(FW)},C(G)).
$$
From this we get an injection
$$
\psi : W\hookrightarrow (FW)^{max},
$$
which is continuous.
Let $\zeta : W^\infty \ra V$ be a morphism in $\CHC_{alg}(G)$.
We get
$$
\alpha(\zeta) : W\hookrightarrow (FW)^{max} \begin{array}{c} \zeta^{max}\\ \longrightarrow\\ {}\end{array} V^{max},
$$
with $\alpha(\zeta)|_{W^\infty}=\zeta$, i.e. the desired surjectivity.
\qed

The Proposition especially implies that any completion of $V\in\CHC_{alg}(G)$ injects into $V^{max}$.
This assertion is immediate for irreducible modules and follows generally by induction on the length.
Thus the use of the term 'maximal completion' is justified.

Let now $\Ga$ denote a cocompact torsion-free discrete subgroup of $G$, let $(\omega ,V_\omega)$ be a finite dimensional unitary representation of $\Ga$ and let $C(\Ga \bs G,\omega)$ denote the space of all continuous functions $f$ from $G$ to $V_\omega$ satisfying $f(\ga x)=\omega(\ga)f(x)$.

For $U,V\in Alg(G)$ and $p,q\ge 0$ let $H_d^q(G,V)$ and $Ext_{G,d}^q(U,V)$ denote the differentiable cohomology and Ext-groups as in \cite{Borel-Wallach}, chap. X.

The next theorem is a generalization of the classical duality theorem \cite{gelf}.

\begin{theorem}\label{higher_duality}
(Higher Duality Theorem) 
Let $\Ga$ be a cocompact torsion-free discrete subgroup of $G$ and $(\omega ,V_\omega)$ a finite dimensional unitary representation of $\Ga$, then for any $V\in \CHC_{alg}(G)$:
\begin{eqnarray*}
H^q(\Ga ,V^{max}\otimes \omega) &\cong& H_d^q(G,Hom_\C(C(\Ga\bs G,\tilde{\omega})^\infty ,V)^\infty)\\
	&\cong& Ext_{G,d}^q(C(\Ga\bs G,\tilde{\omega})^\infty ,V).
\end{eqnarray*}
\end{theorem}

\prf
Let 
$$
\hat{M} := Hom_\C(\tilde{V},C(G)\otimes V_\omega).
$$
Consider the action of $G$ on $\tilde{M}$:
\begin{eqnarray*}
l(g)\alpha (v^*)(x) &:=& \alpha(g^{-1}.v^*)(xg),\\
\end{eqnarray*}

It is clear that
$$
\hat{M}^{l(G)} \cong V^{max}\otimes V_\omega.
$$
Let $M$ denote the algebraic part of $\hat{M}$ with respect to $l$. i.e.
\begin{eqnarray*}
M &=& Hom_\C(\tilde{V},C(G)\otimes V_\omega)^\infty\\
&=& \lim_{\begin{array}{c}\ra\\ K\end{array}} Hom_K(\tilde{V},C(G)\otimes V_\omega)
\end{eqnarray*}
where $K$ runs over the set of all compact open subgroups of $G$.

The group $\Ga$ acts on $M$ by
$$
\ga \alpha(v^*)(x) := \omega(\ga) \alpha(v^*)(\ga x).
$$
This action commutes with $l(G)$.
Now let $\CC^\infty(G)$ be the subcategory of $\CC(G)$ consisting of {\bf smooth} representations.
This means that $V\in\CC(G)$ lies in $\CC^\infty(G)$, when $V^\infty =V$.
The inductive limit topology makes $M$ an element of $\CC^\infty(G\times \Ga)$.
Let $Ab$ denote the category of abelian groups.
The functor $\CC^\infty(G\times \Ga)\ra Ab$ defined by taking $G\times \Ga$-invariants can be written as the composition of the functors $H^0(\Ga,.)$ and $H^0(G,.)$ in two different ways:
$$
\begin{array}{ccc}
\CC^\infty(G\times \Ga) & \begin{array}{c}{H^0(\Ga ,.)}\\ \longrightarrow\\ {}\end{array}& \CC^\infty(G)\\
{ H^0(G,.)}\downarrow& {}& \downarrow H^0(G,.)\\
\CC^\infty(\Ga) &\begin{array}{c}{H^0(\Ga ,.)}\\ \longrightarrow\\ {}\end{array}& Ab,
\end{array}
$$
giving rise to two spectral sequences:
\begin{eqnarray*}
^1E_2^{p,q} &=& H^p_d(G,H^q(\Ga ,M)),\\
^2E_2^{p,q} &=& H^p(\Ga ,H^q_d(G,M)),
\end{eqnarray*}
both abutting to $H^*(G\times \Ga ,M)$.

\begin{lemma}
The module $M$ is $(G,d)$-acyclic and $\Ga$-acyclic.
\end{lemma}

When saying $(G,d)$-acyclic we mean differentiable or strong acyclicity \cite{Borel-Wallach}.
For the discrete group $\Ga$ this notion coincides with usual acyclicity.

\prf
We show the $G$-acyclicity first.
We have
\begin{eqnarray*}
H_d^q(G,M) &=& H_d^q(G,Hom(\tilde{V},C(G)\otimes \omega)^\infty)\\
	&=& Ext_{G,d}^q(\tilde{V},C(G)^\infty\otimes\omega).
\end{eqnarray*}
By \cite{Borel-Wallach} X 1.5 we know that $C(G)^\infty$ is s-injective.
Therefore $H_d^q(G,M)=0$ for $q>0$ as desired.

For the $\Ga$-acyclicity we consider the standard resolution $C^q:=\{ f:\Ga^{q+1}\ra M\}$ with the differential $d: C^q\ra C^{q+1}$ given by
$$
df(\ga_0 ,\dots ,\ga_{q+1}) = \sum_{j=0}^{q+1} (-1)^j f(\ga_0 ,\dots ,\hat{\ga_j}, \dots ,\ga_{q+1}).
$$
The group $\Ga$ acts on $C^q$ by
$$
(\ga f) (\ga_0 ,\dots ,\ga_q) = \ga .f(\ga^{-1}\ga_0 ,\dots ,\ga^{-1}\ga_q).
$$
Then $H^*(\Ga ,M)$ is the cohomology of the complex of $\Ga$-invariants $(C^q)^\Ga$.

Choose an open set $U\subset G$ such that $(\ga U)_{\ga\in\Ga}$ is a locally finite open covering of $G$.
Then there exists a locally constant partition of unity $(\rho_\ga)_{\ga\in\Ga}$ such that $\supp \rho_\ga\subset \ga U$ and $L_{\ga'}\rho_\ga =\rho_{\ga'\ga}$.
Consider a cocycle $f\in(C^q)^\Ga$, $df=0$, $q\ge 1$.
Let $F\in C^{q-1}$ be defined by
$$
F(\ga_0,\dots ,\ga_{q-1})(x) := \sum_{\ga\in\Ga} f(\ga_0 ,\dots ,\ga_{q-1},\ga)(x)\rho_\ga(x),\ \ \ x\in G.
$$
The sum is locally finite in $x$.
A computation shows that $F$ is $\Ga$-invariant and that $dF=(-1)^qf$.
The lemma follows.
\qed

By the lemma both our spectral sequences degenerate.
Since they both abut to $H^*(G\times \Ga ,M)$ we get
\begin{eqnarray*}
H^p(\Ga ,H^0(G,M)) &=& ^1E_2^{p,0}\\
	&=& ^2E_2^{p,0}\\
	&=& H^p_d(G,H^0(\Ga ,M)).
\end{eqnarray*}
The first term is $H^p(\Ga,V^{max}\otimes V_\omega)$.
The last is
$$
H_d^q(G ,Hom_\C (\tilde{V},C(\Ga\bs G,\omega))^\infty).
$$
Since $V$ is of finite length we get an isomorphism as $G$-modules:
$$
Hom_\C(\tilde{V} ,C(\Ga \bs G,\omega))^\infty \cong \bigoplus_{\pi\in\hat{G}} N_{\Ga ,\omega}(\pi) Hom_\C (\tilde{V},\pi^\infty).
$$

We now employ the following

\begin{lemma}
Let $V,W\in\CHC_{alg}$, then we have an isomorphism of $G$-modules:
$$
Hom_\C(V,W)^\infty \cong Hom_\C(\tilde{W},\tilde{V})^\infty.
$$
\end{lemma}

\prf
We have a natural map $Hom_\C(V,W)^\infty \ra Hom_\C(\tilde{W},\tilde{V})^\infty$ given by dualizing: $\alpha\mapsto \tilde{\alpha}$.
This map is injective.
Iterating, we get a map $Hom_\C(\tilde{V},\tilde{W})^\infty \ra Hom_\C(\tilde{\tilde{W}},\tilde{\tilde{V}})^\infty\cong Hom_\C(V,W)^\infty$, and these maps are inverse to each other.
\qed

Recall that $C(\Ga\bs G,\omega)^\infty$ is a direct sum of irreducibles, each occurring with finite multiplicity.
So every $G$-morphism joining $C(\Ga\bs G,\omega)^\infty$ with a module of finite length will factor over a finite sum of irreducibles in $C(\Ga\bs G,\omega)^\infty$.
Therefore the above applies to give
$$
Hom_\C(\tilde{V} ,C(\Ga \bs G,\omega)^\infty) \cong Hom_\C(C(\Ga\bs G,\tilde{\omega})^\infty ,V),
$$
where we additionally have used that the invariant integral puts $C(\Ga\bs G,\omega)^\infty$ in duality with $C(\Ga\bs G,\tilde{\omega})^\infty$.
The theorem is proven.
\qed

Let by a slight abuse of notation $\pi_{\sigma ,\chi ,s_0}$ be the representation parabolically induced from $\sigma \chi |.|^{s_0}$, i.e. the space of this representation consists of measurable functions
$f:G\ra V_\sigma$ such that $f(manx)=\chi(a) |a|^{s_0 +|\rho |}\sigma(m) f(x)$, $f$ is square-integrable on $K$ modulo functions vanishing on a set of measure zero.

The next theorem is the Patterson conjecture as adapted to our situation.

\begin{theorem} 
Assume that the Weyl-group $W(G,A)$ is non-trivial.
The vanishing order of the function $Z_{P,\sigma ,\chi ,\omega}(q^{-s-|\rho|})$ at $s=s_0$ is
$$
(-1)^{q(G)+1}\chi_{_1}(\Ga ,\pi_{\sigma ,\chi ,s_0}^{max} \otimes \omega)
$$
if $s_0\ne 0$ and
$$
(-1)^{q(G)+1}\chi_{_1}(\Ga ,\hat{\pi}_{\sigma ,\chi ,0}^{max} \otimes \omega)
$$
if $s_0=0$, where $\hat{\pi}_{\sigma ,\chi ,0}^{max}$ is a certain extension of $\pi_{\sigma ,\chi ,0}^{max}$ with itself.

Further we have that the usual Euler-characteristic $\chi(\Ga ,\pi_{\sigma ,\chi ,-s_0}^{max} \otimes \omega)$ always vanishes.
\end{theorem}

\prf
So far we know that the vanishing order of $Z_{P,\sigma ,\chi ,\omega}(q^{-s-2|\rho|})$ at $s=s_0$ equals
$$
(-1)^{q(G)+1}\sum_{\pi\in\hat{G}} N_{\Ga ,\omega}(\pi) \sum_{i=0}^{\dim M} (-1)^i \dim H_d^i(M,\pi_{\bar{N}}\otimes\sigma_\chi)_{-s_0}.
$$

Note that $H_d^i(M,\pi_{\bar{N}}\otimes\sigma_\chi)_{-s_0} = H_d^i(M,(\pi_{\bar{N}}\otimes\sigma_\chi)_{-s_0})$.

\begin{lemma}
Let $\pi\in\hat{G}$.
If $s\ne |\rho |$, then $A$ acts semisimply on $(\pi_{\bar{N}})_{s}$.
In any case the $AM$-module $\pi_{\bar{N}}$ has length at most $2$.
\end{lemma}

\prf
Suppose $\pi_{\bar{N}}\ne 0$. Then there is an induced representation $I=Ind_{\bar{P}}^G(\xi\otimes\mu\otimes 1)$ such that $\pi$ injects into $I$.
By the exactness of the Jacquet functor it follows that $\pi_{\bar{N}}$ injects into $I_{\bar{N}}$.
It therefore suffices to show the assertion for induced representations of the form $\pi =I$.
These are given on spaces of locally constant functions $f:G\ra V_{\xi ,\mu}$ such that $f(ma\bar{n}x)=|a|^\mu \xi (m) f(x)$.
The map $\psi : I\ra V_{\xi ,\mu}$, given by $f\mapsto f(1)$ is a $\bar{P}$-homomorphism.
The group $A$ acts on the image by $|a|^\mu$ and this defines an irreducible quotient of $\pi_{\bar{N}}$.
We will show that $\ker \psi$ gives an irreducible subrepresentation of $\pi_{\bar{N}}$ on which $A$ acts by $|a|^{-\mu +2|\rho|}$.
Bruhat-Tits theory tells us that $G$ splits into two $\bar{P}$-double cosets $G=\bar{P} \cup \bar{P}w\bar{P} = \bar{P} \cup \bar{P}w\bar{N}$, where $w$ is the nontrivial element of the Weyl-group $W(G,A)$.
Let $f\in\ker\psi$.
Then $f$ vanishes in a neighborhood of $1$ and thus the integral $\int_{\bar{N}}f(w\bar{n})d\bar{n}$ converges.
(For this recall that by the Bruhat-decomposition $\bar{P}\bs G$ is the one point compactification  of $\bar{N}$.
Therefore $f(1)=0$ implies that the function $\bar{n}\mapsto f(w\bar{n})$ has compact support.) 
This integral defines a map $(\ker\psi)_{\bar{N}}\ra V_{\xi}$.
The group $A$ acts as
\begin{eqnarray*}
\int_{\bar{N}}f(w\bar{n}a)d\bar{n} &=& |a|^{2|\rho|} \int_{\bar{N}}f(wa\bar{n})d\bar{n}\\
	&=& |a|^{2|\rho|} \int_{\bar{N}}f(a^{-1}w\bar{n})d\bar{n}\\
	&=& |a|^{-\mu +2|\rho|} \int_{\bar{N}}f(w\bar{n})d\bar{n}.
\end{eqnarray*}
The claim is now clear.
\qed

To show the theorem consider the case $s\ne -|\rho|$.
Then by semisimplicity it follows
$$
H_d^q(M,(\pi_{\bar{N}}\otimes\sigma_\chi)_{-s_0}) = H^0(A,H_d^q(M,\pi_{\bar{N}}\otimes\sigma_{\chi s_0})),
$$
where $\sigma_{\chi s_0}$ is the $AM$-module given by $\sigma_{\chi s_0}(am)=|a|^{s_0}\chi(a)\sigma(m)$.
Similar to \ref{chichi1} we get
$$
\sum_{q=0}^{\dim M} (-1)^q \dim H^0(A,H_d^q(M,\pi_{\bar{N}}\otimes\sigma_{\chi s_0})) = \sum_{q=0}^{\dim AM} q(-1)^{q+1} \dim H_d^q(AM,\pi_{\bar{N}}\otimes\sigma_{\chi s_0}).
$$
By \cite{Borel-Wallach} p.262 we know
\begin{eqnarray*}
H_d^q(AM ,\pi_{\bar{N}}\otimes \sigma_{\chi s_0})&\cong& H_d^q(AM ,Hom_\C(\tilde{(\pi_{\bar{N}})}, \sigma_{\chi s_0}))\\
	&\cong& Ext_{AM,d}^q((\tilde{\pi})_N,\sigma_{\chi s_0})
\end{eqnarray*}

Now recall that the two Jacquet functors $(.)_N$ and $Ind_P^G(.\otimes 1)$ are adjoint to each other, i.e.
$$
Hom_{AM,d}((.)_N,.)\cong Hom_{G,d}(.,Ind_P^G(.\otimes 1)).
$$
This especially implies that $Ind_P^G(.\otimes 1)$ maps injectives to injectives and that $(.)_N$ maps projectives to projectives.
Further they are both exact, hence preserve resolutions.
It follows that for any $q$:
\begin{eqnarray*}
Ext_{AM,d}^q((\tilde{\pi})_N,\sigma_{\chi s_0})
	&\cong& Ext_{G,d}^q(\tilde{\pi},Ind_P^G(\sigma_{\chi s_0}))\\
	&\cong& Ext_{G,d}^q(\tilde{\pi},\pi_{\sigma ,\chi ,{s_0}-|\rho|}).
\end{eqnarray*}

Therefore we have for the order
\begin{eqnarray*}
ord_{s=s_0} Z_{P,\sigma ,\chi ,\omega}(q^{-s-2|\rho|})
&=& (-1)^{q(G)+1}\sum_{\pi\in\hat{G}} N_{\ga ,\omega}(\pi) \sum_{q=0}^{\dim AM} q(-1)^{q+1} \dim Ext_{G,d}^q(\tilde{\pi},\pi_{\sigma ,\chi ,s_0-|\rho|})\\
	&=& (-1)^{q(G)+1}\sum_{q=0}^{\dim AM} q(-1)^{q+1} \dim Ext_{G,d}^q(C(\Ga\bs G,\tilde{\omega})^\infty,\pi_{\sigma ,\chi ,s_0-|\rho|})\\
&=& (-1)^{q(G)+1}\sum_{q=0}^{\dim AM} q(-1)^{q+1} \dim H^q(\Ga,\pi_{\sigma ,\chi ,s_0-|\rho|}^{max}\otimes\omega),
\end{eqnarray*}
the last equality comes by \ref{higher_duality}.

In the case $s=-|\rho|$ finally, one replaces $A$ by $A^2=v^{-1}(2\Z)$, the set of all $a\in A$ with even valuation.
Then $\pi_{\sigma ,\chi ,s_0-|\rho |}$ is replaced by $Ind_{A^2MN}^G(\sigma_{\chi s_0})$, which is the desired extension.
\qed

The method employed here for the $p$-adic case also gives a new proof of the Patterson conjecture in the  setting of real Lie groups \cite{BuOl}.

$ $

{\small  Math. Inst. d. Univ., INF 288, 69126 Heidelberg, Germany; e-mail: anton@mathi.uni-heidelberg.de}
\end{document}